\documentclass[onefignum,onetabnum]{siamart190516}

\usepackage{lipsum}
\usepackage{amsfonts}
\usepackage{graphicx}
\usepackage{epstopdf}
\usepackage{algorithmic}
\usepackage{amsopn}
\usepackage{fancybox}
\usepackage{mathtools}
\usepackage{todonotes}
\usepackage{amssymb}
\usepackage{pgf,tikz}
\usepackage{pgfplots}
\usepackage{color}
\pgfplotsset{compat=newest}
\usepackage{mathabx}
\usepackage{hyperref}
\usepackage{cleveref}

\newsiamremark{remark}{Remark}
\newsiamremark{hypothesis}{Hypothesis}
\crefname{hypothesis}{Hypothesis}{Hypotheses}
\newsiamthm{claim}{Claim}

\headers{QN Optimization with Hessian Samples}{J.~Azzam, D.~Henderson, B.~W.~Ong and A.~A.~Struthers}

\title{QN Optimization with Hessian Samples}

\author{ Joy~Azzam \thanks{Michigan Technological University, Houghton, MI
  (\email{atazzam@mtu.edu}).}
  \and 
  Daniel P.~Henderson \thanks{Michigan Technological University, Houghton, MI 
  (\email{dphender@mtu.edu}).}
  \and
  Benjamin W.~Ong \thanks{Michigan Technological University, Houghton, MI
    (\email{ongbw@mtu.edu}), \url{http://mathgeek.us/}.}
  \and
  Allan A.~Struthers \thanks{Michigan Technological University, Houghton, MI
  (\email{struther@mtu.edu}).}
}

\DeclarePairedDelimiterX{\Finner}[2]{\langle}{\rangle_F}{#1, #2}

\DeclareMathOperator*{\argmin}{arg\,min}

\ifpdf
\hypersetup{
  pdftitle={QN Optimization with Hessian Samples},
  pdfauthor={J.~Azzan, D.~Henderson, B.~W.~Ong and A.~A.~Struthers}
}
\fi

\begin{document}
\newcommand{\gAD}{\mbox{\tt gAD}}
\newcommand{\gHS}{\mbox{\tt gHS}}
\newcommand{\QN}{\mbox{\tt QN}}
\newcommand{\TR}{\mbox{\tt TR}}

\maketitle

\begin{abstract}
This article explores how to effectively incorporate 
curvature information generated using SIMD-parallel 
forward-mode Algorithmic Differentiation (AD) into 
unconstrained Quasi-Newton (QN) minimization of a smooth objective function,
$f$. Specifically, forward-mode AD can be used to generate block Hessian samples 
$Y=\nabla^2 f(x)\,S$ whenever the gradient is evaluated. Block QN algorithms then 
update approximate inverse Hessians, $H_k \approx \nabla^2 f(x_k)$, with
these Hessian samples. Whereas standard line-search based BFGS algorithms 
carefully filter and correct secant-based approximate curvature information
to maintain positive definite approximations, our algorithms 
directly incorporate  Hessian samples to update indefinite
inverse Hessian approximations without filtering. The sampled 
directions supplement the standard QN two-dimensional trust-region sub-problem
to generate a moderate dimensional subproblem which can exploit negative curvature.
The resulting quadratically-constrained quadratic program is solved 
accurately with a generalized eigenvalue algorithm and the step 
advanced using standard trust region step acceptance and radius adjustments.
The article aims to avoid serial bottlenecks, exploit accurate positive
and negative curvature information, and conduct a preliminary evaluation 
of selection strategies for $S$.

\end{abstract}

\begin{keywords}
  Optimization,
  Randomized algorithms,
  Samples,
  Quasi-Newton
\end{keywords}

\begin{AMS}
  68W20, 
  68W2, 
  65F35, 
  90C53 
\end{AMS}

\section{Literature Discussion}
Schnable \cite{SchnabelConcurrentFnEvals87} first discusses incorporating 
parallel function evaluations to improve Hessian approximation
in optimization algorithms. Together with his colleagues, Byrd, Schnable, and Shultz 
\cite{ByrdSchnableShultzParallelFnEvals88,ByrdSchnableShultzParallelQN88}
supplement several standard Quasi Newton optimization schemes with a small 
number of finite difference second derivative approximations and conclude that
the supplemental information improves all the variants considered.  Their study includes
the underlying QN update, the update order within each step,  and 
a variety of strategies to choose a few supplemental directions. In their conclusions
Byrd, Schnable, and Shultz \cite{ByrdSchnableShultzParallelQN88} suggest 
supplementing the 
extremely simple Symmetric-Rank-One (SR1) QN update with additional second derivative
information: the motivation for this suggestion is that Conn, Gould and Toint 
find in \cite{ConnGouldTointSR1Good88} that SR1 is better than 
Powell-Symmetric-Broyden (PSB), 
Davidon-Fletcher-Powell (DFP), and 
Broyden–Fletcher–Goldfarb–Shanno (BFGS). 
The review article by Schnable \cite{SCHNABEL_Rev_ParComp_1995}
summarizes this and other early parallel optimization approaches. 
More recently, Gau and Goldfarb \cite{BLOCKBFGSGoldfarb_2018} 
implemented and tested line-search based algorithms using block BFGS 
updates on subsets of previous directions, and a family of Quasi-Newton 
algorithm \cite{QNGoldfarb_2018_NoLineSearch} which avoids a line-search 
for a restricted class of cost functions. 
In yet another approach, Berahasa, Jahanib, Richtarik and Takac \cite{QN_NoMemory} 
develop zero-memory Block-BFGS and block-SR1 algorithm using
only one (the most recent)  set of AD generated Hessian samples.
Their BFGS variant is implemented with a Line-Search 
(requiring a positive definite approximate Hessian) while 
their SR1 variant is implemented with a Conjugate Gradient based 
approximate solution of a full dimensional trust region sub problem. 

The review article \cite{NonLinOptParComReview2003} and the articles in 
the associated special edition of {\em Parallel Computing} emphasize 
improving the parallel efficiency of the underlying linear algebra and 
introducing parallelism through simultaneous 
(with different starting points and/or different but possibly related 
QN updates) 
line searches. The second-order section of the extensive review article 
\cite{SiamRev2018NocedalOptMachLearn} provides a more recent update with a
focus on algorithms for very high dimensional problems. 

A number of recent developments make it appropriate to revisit the topics in 
\cite{ByrdSchnableShultzParallelQN88} with modern computational tools.  
The GPU-enabled forward-mode  Algorithmic Differentiation (AD), implemented in 
the open-source computational tool Julia \cite{JuliaAD2016} and other software projects, 
can efficiently replace the finite difference approximations used 
in \cite{ByrdSchnableShultzParallelQN88} and greatly increases
the number of simultaneous Hessian samples. The generalized 
eigenvalue based trust-region sub-problem solver 
developed by Adachi et al.~\cite{GEVSubProblem2017} 
can replace the line search with an accurate and robust 
moderate-dimensional trust-region sub-problem (TRSP), implemented in   
Julia \cite{GEVActiveSet21}. 

The goal of this manuscript is to incorporate SIMD-parallel Hessian samples into
QN updates to generate provably convergent QN like algorithms. The proposed 
algorithms avoid serial bottlenecks by using indefinite approximate Hessians.
The standard two-dimensional TRSP minimizes a quadratic model over the span
of the steepest-descent and Newton directions. This 2D search space is extended 
(through a specific Hessian re-sampling strategy) to include additional supplementary 
directions with accurate curvature information on the resulting moderate dimensional 
sub-space. With standard trust-region controls this gives a provably convergent algorithm
with rapid asymptotic convergence to non-degenerate local minimizers. 

The article explores two simple indefinite updates (a block variant of
SR1 and a block variant of Powell Symmetric Broyden) which can be
directly implemented on accurate AD curvature information. 
In contrast, line-search based block methods (such as a block BFGS
or DFP) need to carefully filter and correct approximate curvature information to maintain 
positive definite Hessian approximations.
The article explores how the selection strategy for and number of 
supplementary directions affects the algorithms based on these two standard
indefinite QN updates. 

\Cref{sec:Notation} introduces essential assumptions and notation.
\Cref{sec:AD}
discusses the choices made to evaluate curvature information
using Algorithmic Differentiation (AD). 
\Cref{sec:QN} discusses the advantages of 
indefinite QN updates for trust-region optimization,
explains why block BFGS and DFP are not suitable, 
and presents the simple and highly-parallel block SR1 and PSB updates used. 
\Cref{sec:TRSP} describes the trust-region sub-problem underlying the
algorithm. 
\Cref{sec:SUpdate} explains how new supplemental directions are chosen.
\Cref{sec:Alg} describes the assembled algorithm including: 
a  simple mean curvature estimate used to initialize $H$; 
a simple initial trust region radius $\Delta_0$ based on the curvature 
in the steepest descent direction;
standard trust region control; cost evaluation; and pseudo code. 
\Cref{sec:Experiments} describes our numerical experiments.
\Cref{sec:Conclusions} summarizes the results and
future plans. 

\section{Notation and Assumptions} 
\label{sec:Notation}
We assume throughout we are seeking the unconstrained minimum of a $C^2$ 
function $f: \mathbb{R}^n \rightarrow \mathbb{R}$, with
gradient
$g(x) = \nabla f(x)$
given by $g: \mathbb{R}^n \rightarrow \mathbb{R}^n$, 
and that we can efficiently sample the Jacobian of $g$ 
(which is the Hessian of $f$) 
$J(x) = \nabla g(x) = \nabla^2 f(x)$ 
using forward-mode AD.  In this context, sample the Hessian means that 
whenever we compute $g(x)=\nabla f(x)$ we can efficiently and 
simultaneously compute $Y = J(x) \, S \in \mathbb{R}^{n\times w}$
for a block of $w$ directions, $S \in \mathbb{R}^{n\times w}$.
We write {\tt orth}($M$) for an orthogonalization of $M$ 
(implemented as {\tt Matrix(qr(M).Q)} in Julia),
$M^\dagger_\delta$
(implemented as {\tt pinv}($M$; rtol = $\delta$) in Julia)
for the $\delta $ thresholded pseudo-inverse of $M$,
and $M \sim \mathcal{N}_{0,1}^{n \times w}$ 
(implemented as {\tt randn}(n,w))
for an $n \times w$ matrix with elements drawn from the standard 
normal distribution. 
We denote the current search point 
by $x_k$ with:  
objective function value $f_k =f(x_k)$; 
gradient $g_k = g(x_k) = \nabla f(x_k)$;
Jacobian $J_k = J(x_k) = \nabla^2 f(x_k)$;
and $h_k = h(x_k)$
where $h(x) = J(x) \, g(x) = \nabla^2 f(x) \nabla f(x)$.
Symmetric Quasi-Newton (QN) approximations $B_k$ and $H_k$ 
are updated (using indefinite updates which can incorporate 
negative curvature information) to satisfy $B_k \approx J_k$ and 
$H_k \approx J_k^{-1}$. 
The Frobenius inner product of matrices $A$ and $B$
is denoted by $\left<A,B\right>_F$ and the Frobenius norm 
$||A||_F^2 = \left<A, A\right>_F$ is used throughout.

\section{Algorithmic Differentiation and Hessian Samples} 
\label{sec:AD}
In our context, sampling the Hessian means that 
whenever we compute $g(x)=\nabla f(x)$, we can efficiently and 
simultaneously compute $Y = \nabla^2 f(x) \, S \in \mathbb{R}^{n\times w}$
for a block of $w$ directions $S \in \mathbb{R}^{n\times w}$.
The Julia {\tt ForwardDiff} package 
\cite{JuliaAD2016} 
modifies the code of $g$ to code 
$\gAD (x,S): \mathbb{R}^n \times\mathbb{R}^{n\times w} \rightarrow 
\mathbb{R}^n \times \mathbb{R}^{n\times w}$
for the simultaneous combined computation.
In practice, $\gAD$ is embarrassingly SIMD parallel. 
Provided $w$ does not exceed the available processor resources, 
evaluating $\gAD(x,S)$ takes only $2-3$ times as long as evaluating $g(x)$. 
Since modern GPUs have over $256$ cores organized into SIMD warps of 
8, 16 or 32 threads, values of $w\leq 256$ are feasible on most 
commodity modern hardware with much larger values feasible on specialized hardware. 

The algorithm we will generate centers around generating 
accurate curvature information at a single point
for the steepest descent direction, in addition to other sampled directions.
To accomplish this, the curvature information 
Incorporated in \cref{alg:gHS} is generated by 
two sequential $\gAD$ calls. A first call to 
$\gAD$ computes the gradient and a first
Hessian sample $(g, Y_1) \leftarrow \gAD(x,S_1)$.  
The gradient, $g$, is included in a second direction set 
$S_2$ so that the Hessian sample $Y_2$ computed by 
$\gAD(x,S_2)$ contains $h(x) = J(x)\, g(x) = \nabla^2 f(x) \,\nabla f(x)$. 
This can be organized in a number of ways. \Cref{alg:gHS} gives 
the simple (and almost certainly non-optimal)
choices made for a combined gradient and 
Hessian Sample operation $\gHS(x,S)$.

\begin{algorithm}[ht]
    \caption{Gradient and Hessian Sample: 
    $(g, h, Y) \leftarrow \gHS(x, S)$
    }
    \begin{algorithmic}[1]
    \REQUIRE{
    $x \in\mathbb{R}^{n}$ and $S \in\mathbb{R}^{n \times (2w-1)}$.}
    \STATE Compute $\left( g, Y_1 \right) \leftarrow \gAD(x,S\left[\; : \;,\; 1:w\;\right])$
    \COMMENT{Input the first $w$ columns of $S$ to $\gAD$}
    \STATE Compute $\displaystyle \left( g, Y_2 \right) 
    \leftarrow \gAD\left(x,\left[\;S\left[\;: \;,\; w+1\!:\!\text{end}\;\right]\,,\,g\;\right]\right)$
    \COMMENT{Input the last $w-1$ columns of $S$ and $g$ to \gAD} 
    \STATE{Assemble $Y \leftarrow [\;Y_1 \;,\; Y_2[\;: \;, \; 1:\text{end}-1]\;]$}
    \RETURN $\left( g, Y_2[\;:\;,\;\text{end}\;] , \; Y \right)$ 
    \COMMENT{By construction, $Y_2[\;:\;,\;\text{end}\;] = \nabla^2 f(x)\, g$}
    \end{algorithmic}
    \label{alg:gHS}
\end{algorithm}

\section{Quasi-Newton Updates}
\label{sec:QN}
Several Quasi-Newton updates (with a variety of update details)
are tested in \cite{ByrdSchnableShultzParallelFnEvals88}. 
The subsequent article \cite{ByrdSchnableShultzParallelQN88}
which focuses on the now dominant BFGS algorithm 
notes that 
an indefinite SR1 update might work well,
as evidenced in an almost contemporaneous article
\cite{ConnGouldTointSR1Good88}.
In this preliminary study we use the 
indefinite block SR1 and PSB algorithms to 
incorporate a block of accurate Hessian samples 
$Y = \nabla^2 f(x_{k+1}) \, S$ after successful steps.
We explore the effect of including a standard secant 
curvature estimates 
$\left( \nabla f(x_{k+1}) - \nabla f(x_k)\right) 
\approx \nabla f^2 (x_k) \left( x_{k+1}- x_k \right)$ 
before the block update and including the prior step in
the block update.
The Hessian is not updated on a failed step 
since the gradient is not evaluated and there
is no new curvature information.

Conn, Gould and Toint 
\cite{ConnGouldTointSR1Good88} present evidence
that when the underlying Hessian is indefinite, SR1 
(with suitable globalization strategies for indefinite 
Hessian approximations) can outperform updates (like BFGS) 
designed to maintain positive definite Hessian approximations.
We use a globalization strategy (the moderate dimensional
trust region presented in \Cref{sec:TRSP}) which can
exploit negative curvature and select 
QN updates suitable for indefinite approximations. 

Block versions of updates which implicitly use the 
current Hessian (such as BFGS and DFP) require care 
when the Hessian is indefinite.  
We explain the issues for the block BFGS update 
in \cref{alg:BFGS} 
\cite{ByrdSchnableShultzParallelQN88,ByrdSchnableShultzParallelFnEvals88}
which incorporates a block of curvature estimates 
$\widehat{V}\approx \nabla^2 f(x) \, U$
into an inverse Hessian approximation $H^{-1} \approx \nabla^2 f(x)$.
\begin{algorithm}[ht]
    \caption{block BFGS Update: $H_{\text{BFGS}} \leftarrow {\tt BFGS}(H,U,\widehat{V},\delta)$.}
    \begin{algorithmic}[1]
    \label{alg:BFGS}
    \REQUIRE{
    SPD $H \in\mathbb{R}^{n\times n}$ and
    $U,\widehat{V}\in\mathbb{R}^{n \times 2w}$ with 
    $U^\top \widehat{V}$ approximately symmetric}
    \STATE{
    Filter and correct $\widehat{V}$ 
    to $V$ consistent with $V =A \,U$ for some 
    SPD $A$.}
    \STATE{Compute $T \leftarrow \left( U^\top V \right)^\dagger_\delta$.} 
    \RETURN{$U \, T \, U^\top +
    \left(I- U \, T \, V^\top \right) H 
    \left(I- V \, T \, U^\top \right)$. }
    \end{algorithmic}
\end{algorithm}
There are two 
reasons for the initial filtering step in line  of \cref{alg:BFGS}.
Finite difference approximations and/or multiple 
secant estimates generate approximations 
$\widehat{V} \approx \nabla^2 f(x) \,U$, however,
 $\widehat{V}$ needs to be corrected to ensure 
$U^\top \widehat{V}$ is symmetric.
When $\nabla^2 f(x)$ is not positive definite,
negative curvature directions in $\widehat{V} \approx \nabla^2 f(x) \,U$
need to be filtered to ensure $H_{\text{BFGS}}$ is SPD. 
These corrections and filters
\cite{ByrdSchnableShultzParallelQN88,ByrdSchnableShultzParallelFnEvals88,QN_NoMemory}
are inherently serial. In 
contrast, there is no need to correct accurate AD Hessian 
samples $V = \nabla^2 f(x) \, U = J(x)\, U$ from $\gHS(x,S)$ or
filter negative curvature directions for trust region 
based algorithms which can use indefinite Hessian approximations
$H$. 
Note the BFGS update is a critical point of
\begin{equation}
\label{eq:BFGS-OPT}
    \argmin_{A \,V =U, A=A^\top} \left<\; (H - A) \,\mathcal{J}^{-1} \;,\; (H-A) \; \right>_F
\end{equation}
for any invertible $\mathcal{J}$ consistent with the sample in the sense
that $U^\top V = U^\top \mathcal{J} \, U$.  
Running BFGS without filtering produces meaningless
updates since if $U^\top V$ is not SPD \cref{eq:BFGS-OPT} gives a saddle 
point of an unbounded minimization problem.

We test two indefinite block QN updates with identical trust-region sub-solvers and 
controls. \Cref{alg:SR1} specifies block Symmetric Rank 1 (block SR1) which is a direct block generalization 
\cite{QN_NoMemory} of the rank one SR1 update: block SR1 is the algebraically minimal update $H_{\text{SR1}}$ which 
satisfies the block inverse secant condition $ H_{\text{SR1}} V=U$.  
\Cref{alg:PSB} specifies block Powell-Symmetric-Broyden (block PSB) which is a direct block generalization of 
the rank two Powell-Symmetric-Broyden (PSB) update \cite{ByrdSchnableShultzParallelFnEvals88,ByrdSchnableShultzParallelQN88}:
block PSB is the minimal Frobenius norm change satisfying the 
block inverse secant condition $ H_{\text{PSB}} V=U$, i.e.,
\begin{equation*}
    H_{\text{PSB}} = \argmin_{A V =U, A=A^\top}||H - A||_F
    = \argmin_{A V =U, A=A^\top} \left<\; H - A \;,\; H-A\; \right>_F.
\end{equation*}
\begin{algorithm}[htbp]
    \begin{algorithmic}[1]
    \REQUIRE{
    $H \in\mathbb{R}^{n\times n}$ with $H=H^\top$;
    $U,V \in\mathbb{R}^{n \times 2w}$ with $U^\top V = V^\top U$;
    $\delta > 0$.}
    \STATE{Compute $T\leftarrow((U-H \, V)^\top V)^\dagger_\delta$.} 
    \RETURN{$H + \left(U-H \, V\right) T \left(U-H \,V\right)^\top$. }
    \end{algorithmic}
    \caption{block SR1 Update: $H_{\text{SR1}} \leftarrow $ {\tt SR1}$(H,U,V,\delta)$.}
    \label{alg:SR1}
\end{algorithm}
\begin{algorithm}[ht]
    \caption{block PSB Update: $H_{\text{PSB}} \leftarrow {\tt PSB}(H,U,V,\delta)$.}
    \begin{algorithmic}[1]
    \label{alg:PSB}
    \REQUIRE{
    $H \in\mathbb{R}^{n\times n}$ with $H=H^\top$;
    $U,V \in\mathbb{R}^{n \times 2w}$ with $U^\top V = V^\top U$;
    $\delta > 0$.}
    \STATE{Compute $T_1 \leftarrow \left( V^\top V \right)^\dagger_\delta$.} 
    \STATE{Compute $T_2 \leftarrow V T_1 (U-H V)^\top$.} 
    \RETURN{$H + T_2 +T_2^\top - T_2 \, V \, T_1 \, V^\top $. }
    \end{algorithmic}
\end{algorithm}
\Cref{alg:SR1,alg:PSB} do not need to filter 
accurate curvature information from $\gHS$. The  pseudo-inverse tolerance 
$\delta=10^{-12}$ simply restricts excessively large updates.

Byrd and Schnabel 
\cite{ByrdSchnableShultzParallelQN88, ByrdSchnableShultzParallelFnEvals88}
considered various additional updates and discovered that including the approximate 
secant curvature information
\begin{equation}
\label{eq:SecantUpdate}
    y_{k+1}=(\nabla f_{k+1} - \nabla f_k) \approx \nabla^2 f(x_{k+1}) p_k = 
    \nabla^2 f(x_{k+1})(x_{k+1} - x_k) 
\end{equation}
improved their algorithms. They recommend a preliminary QN update with 
the approximate curvature information
\cref{eq:SecantUpdate}
before a second update to incorporate the additional curvature information.
We perform numerical experiments which replicate this observation and
evaluate a possible block replacement.

\section{Trust Region Sub-Problem}
\label{sec:TRSP}
Trust region algorithms are based on approximate solutions of the 
$n$ dimensional quadratically constrained quadratic program
\begin{equation}
\label{eq:nDimTRp}
    p_k = \argmin_{|p|\leq \Delta_k}
    \frac{1}{2} p^\top H_k^{-1} \, p + \nabla f_k^\top p
    \quad \mbox{where} \quad H_k^{-1} \approx \nabla^2 f(x_k).
\end{equation}
If $H_k$ is full rank \cref{eq:nDimTRp} is equivalent 
(with $p_k = H_k q_k$) to 
\begin{equation}
\label{eq:nDimTRq}
    q_k = \argmin_{|H_k q|\leq \Delta_k} m_k(q) 
    \quad \mbox{where} \quad
    m_k(q) = \frac{1}{2} q^\top H_k \, q  + (H_k \nabla f_k)^\top q.
\end{equation}
The standard two-dimensional subspace approximation 
(discussed on p76 of \cite{NocedalAndWright06})
minimizes  \cref{eq:nDimTRp} for  
$p = a_1 \nabla f_k + a_2 H_k \nabla f_k$.

The sampling procedure $\gHS(x_k,S_k)$ generates accurate curvature 
information $[\;Y_k\;,\; h_k\;] = \nabla^2 f(x_k) \, [\;S_k\;,\;\nabla f_k\;]$ in the 
directions specified by the columns of $S_k$  and $\nabla f_k$
and the inverse Hessian approximation, $H_k$, is 
immediately updated 
(using either \cref{alg:BFGS} or \cref{alg:SR1})
to match with 
\begin{equation*}
U=\left[\; S_k \;,\; \frac{\nabla f_k}{\|\nabla f_k \|}\; \right]
\quad \mbox{and} \quad
V= \left[\;Y_k \;,\; \frac{h_k}{\|\nabla f_k \|}\; \right]
\end{equation*}
The scaling weights all the columns of $V$ equally and
maintains a well-conditioned computation when 
$\nabla f_k$ is large or small. 
Thus,
\begin{equation}
\label{eq:ExactOnSamples}
    Y_k = H_k^{-1} S_k 
    \quad \mbox{and} \quad
    h_k = H_k^{-1} \nabla f_k. 
\end{equation}
The standard 2D approximation
($p$ in the span of $\nabla f_k$ and $H_k \nabla f_k$)
is expanded with the columns of $S_k$ to give
the explicit representation
$p = M_k \,a$ where 
\begin{equation*}
    M_k = \left[\;\frac{\nabla f_k}{\|\nabla f_k\|} \;,\; \frac{H_k \nabla f_k}{\|\nabla f_k\|} \;,\; S_k\; \right]
    \in \mathbb{R}^{n\times(2w+1)},
    \quad
    a \in \mathbb{R}^{2w+1}.
\end{equation*}
In terms of $q=H_k^{-1}p$, the equivalent representation 
(since $H_k$ is exact on the sample \cref{eq:ExactOnSamples}) 
is $H_k^{-1}M_k$. We use the orthogonal representation
$Q_k = \mbox{\tt orth} \left( H_k^{-1} M_k \right)$
where 
\begin{align}
\label{eq:Qkrep} 
    H_k^{-1} M_k &= \left[\; H_k^{-1}\frac{\nabla f_k}{\|\nabla f_k\|}
    \;,\; 
    \frac{\nabla f_k}{\|\nabla f_k\|} \;,\; H_k^{-1} S_k \;
    \right]
    = \left[ \;
    \frac{h_k}{\|\nabla f_k\|} \;,\; \frac{\nabla f_k}{\|\nabla f_k\|}
    \;,\; 
    Y_k \;\right].
\end{align}
Thus, our new trust-region approximation (with $Q_k$ from \cref{eq:Qkrep} and
$m_k$ from \cref{eq:nDimTRq}) is 
\begin{equation}
      \label{eq:TRSubProb} 
      a_{k} = 
      \argmin_{a\in\mathbb{R}^{2w+1}: |H_k Q_k \,a|\leq \Delta_k} m_k(Q_k \,a)
\end{equation}
giving the trial step $x_{k+1}=x_k + p_k = x_k + H_k \,Q_k \,a_k $. 
The Julia {\tt trs} package \cite{JuliaTRSPackage} (developed for 
\cite{GEVActiveSet21} based on 
\cite{GEVSubProblem2017}) computes accurate eigen-value based 
solutions of
\begin{equation*}
    a_k = \argmin_{a^\top C a \; \leq \; \Delta_k^2}\frac{1}{2} a^\top P a + 
    b^\top a.
\end{equation*} 
We use the robust small-scale solver (based on a dense 
generalized eigenvalue decomposition) and compute accurate
solutions of \cref{eq:TRSubProb} with
\begin{equation}
\label{eq:JuliaTRS}
    a_k = \mbox{{\tt trs\_small}}(P,b,\Delta_k,C)
    \end{equation}
where the arguments are 
\begin{equation*}
    P = Q_k^\top \,H_k\, Q_k,
    \quad
    b =  Q_k^\top\, H_k\, \nabla f_k, 
    \quad
    \mbox{and}
    \quad
    C = Q_k^\top\, H_k^2\, Q_k. 
\end{equation*}

\section{Supplemental Directions}
\label{sec:SUpdate}
Lastly, we need to address how to select the  
supplementary directions, 
$S \in \mathbb{R}^{n \times(2w-1)}$, in \cref{alg:gHS}.
The six supplemental direction 
variants considered are:
\begin{subequations}
    \begin{align}
        \label{eq:SUpdate1}
        S_{k+1} &= {\tt orth}\left( M \right),
        \mbox{ where } 
        M \sim \mathcal{N}_{0,1}^{n \times (2w-1)}; \\
        \label{eq:SUpdate2}
        S_{k+1} &= {\tt orth}\left(M - S_{k} (S_{k}^\top M) \right),
        \mbox{ where } 
        M \sim \mathcal{N}_{0,1}^{n \times (2w-1)};\\
        \label{eq:SUpdate3}
        S_{k+1} &= {\tt orth}\left(Y_{k} - S_{k} (S_{k}^\top Y_{k}) \right);\\
        \label{eq:SUpdate4}
        S_{k+1} &= {\tt orth} \left( \left[ \; {\tt orth}\left( M \right) \;,\; p_k \; \right] \right),
        \mbox{ where } 
        M \sim \mathcal{N}_{0,1}^{n \times (2w-2)}; \\
        \label{eq:SUpdate5}
        S_{k+1} &=  {\tt orth} \left( \left[ \;  {\tt orth}\left(M - S_{k} (S_{k}^\top M) \right) \;,\; p_k \; \right] \right),
        \mbox{ where } 
        M \sim \mathcal{N}_{0,1}^{n \times (2w-2)};\\
        \label{eq:SUpdate6}
        S_{k+1} &=  {\tt orth} \left( \left[ \; 
        {\tt orth}\left(Y_{k}[\;:\;,\;1:\text{end}\!-\!1] - S_{k} (S_{k}^\top Y_{k}[\;:\;,\;1:\text{end}\!-\!1]) \right)
        \;,\; p_k \; \right] \right).
    \end{align}
\end{subequations}
The idea behind \cref{eq:SUpdate2} was to prevent immediate re-sampling 
(which will happen in the simple randomization \cref{eq:SUpdate1})
by orthogonalizing against the immediate previous directions. 
The idea behind \cref{eq:SUpdate3} was to guide the algorithm to 
accurately resolve eigen-space associated with the larger Hessian eigenvalues.
As noted in \cref{sec:QN}, Byrd and Schnabel 
\cite{ByrdSchnableShultzParallelQN88, ByrdSchnableShultzParallelFnEvals88}
perform a preliminary secant update to incorporate the approximate secant curvature information
along the previous step $p_k = x_{k+1}-x_k$ from \cref{eq:SecantUpdate}.
A simple block alternative is to include in the
$p_k$ in $S_{k+1}$ which ensures that the accurate curvature
$\nabla^2 f(x_{k+1})\,p_k$ is incorporated in the inverse Hessian.  
\Cref{eq:SUpdate4,eq:SUpdate5,eq:SUpdate6}
are simply variants of 
\cref{eq:SUpdate1,eq:SUpdate2,eq:SUpdate3}
which include $p_k$.

\section{Algorithmic Overview, Motivation, and Details}
\label{sec:Alg}
The primary goals when designing the algorithm was
to extract maximal benefit from AD generated
curvature information while avoiding linear solves
in the potentially large ambient dimension $n$.  
Secondary goals (which drove many of the details) 
were a clean flow of information and 
provably better objective function reduction 
(at each step) than familiar convergent benchmarks. \Cref{alg:Alg} has pseudo code for the assembled algorithm. Some comments are in order.

Line 6: QN algorithms are commonly initialized with a multiple of the 
identity.  The mean of the eigenvalues of 
$S_0^\top \nabla f^2 (x_0)\, S_0 = S_0^\top Y_0$
where $Y_0 = \gHS(x_0,S_0)$ provides a natural estimate 
for this multiplier.  This initialization
is immediately updated with the QN update to 
give $H_0$ satisfying 
$Y_0 = H_0 S_0$.


Lines 5, 12 \& 13: We adopt the simple standard strategy and terminology
from \cite{NocedalAndWright06} for updating the radius $\Delta_k$
and accepting or rejecting the trial step
$x_k + H_k \,Q_k\, a_k$.
Model quality is assessed by measuring the ratio 
of the actual decrease  of the objective function, ${f_k - f(x_k + H_k\, Q_k\, a_k)}$, 
to the model decrease,  $m_k(0) - m_k(Q_k\, a_k)$, viz.,
\begin{equation}
\label{eq:rho}
     \rho = 
    \frac{f_k - f(x_k + H_k \,Q_k \,a_k) }{m_k(0) - m_k(Q_k \,a_k)}.
\end{equation}
The trust-region radius is updated 
(our experiments use the large maximum trust region radius 
$\Delta_{\mbox{max}}= 100$)
as follows
\begin{equation}
\label{eq:TRUpdate}
      \Delta_{k+1} = \left\{ 
    \begin{array}{rl}
         0.25 \Delta_k 
         & \mbox{if } \rho < 0.25 \\
         \min \left(2 \Delta_k, \Delta_{\mbox{max}}\right)
         & \mbox{if } \rho > 0.75 \mbox{ and } ||H_k Q_k a_k||=\Delta_k \\
         \Delta_k & \mbox{otherwise}
    \end{array}
    \right. .
\end{equation}
We reject the step if $\rho \leq 0$ (this is the 
standard \cite{NocedalAndWright06} trust-region control 
with rejection parameter $\eta=0$): we retain all 
previous values with the updated subscript $k+1$ except 
the radius $\Delta_k$ and recompute \cref{eq:JuliaTRS} 
with the $\Delta_{k+1}=0.25\Delta_k$ since $\rho < 0.25$. 
We accept the step if $\rho > 0$. We set $x_{k+1} = x_k + H_k Q_k a_k$ and 
$f_{k+1} = f(x_{k+1})$ then resample Hessians and update QN approximations. The initial trust-region radius is set to 
$1.1\times$ the distance to the isotropic 
    quadratic critical point with mean curvature $\alpha_0$. 
    In our experiments $\Delta_{\tt max} = 100$.
    
\begin{algorithm}[tbp]
    \caption{Trust Region Quasi-Newton Optimization
    with Hessian Samples}
    \begin{algorithmic}[1]
    \label{alg:Alg}
    \REQUIRE{
    Convergence tolerance $\epsilon > 0$;
    Pseudo-Inverse tolerance $\delta>0$;
    Initial point $x_0\in\mathbb{R}^{n}$;
    function handle $f$;
    Preliminary QN update flag, {\em pflag}.}
    \STATE{Compute $f_0 \leftarrow f(x_0)$; Draw
    $M \sim \mathcal{N}_{0,1}^{n \times (2w-1) }$;
    Set $S_0 \leftarrow {\tt orth} \left(M \right)$}.
    \STATE{Compute 
    $\left( \nabla f_0, h_0, Y_0 \right) \leftarrow \gHS (x_0,S_0)$. }
    \COMMENT{see \Cref{alg:gHS}}
    \STATE{Construct 
    $\displaystyle U_0 = \left[\;S_0 \;,\; \frac{\nabla f_0}{\|\nabla f_0\|} \;\right]$
    and  $\displaystyle V_0 = \left[\; Y_0 \;,\; \frac{h_0}{\|\nabla f_0\|}\;\right]$. }
    \STATE{Compute initial mean curvature estimate 
    $\alpha_0 = {\tt mean}({\tt eig}(S_0^\top Y_0))$.}
    \STATE{Compute initial 
    $\Delta_0 = \min
    \left( 1.1 \times \frac{\|\nabla f_0\|}{2 |\alpha_0|} \;,\;
    \Delta_{\tt max}=100 \right)$.}
    \STATE{Initialize 
    $H_0 \leftarrow \QN \left( \alpha^{-1} I, U_0, V_0,\delta 
    \right)$.}
    \COMMENT{See \Cref{alg:BFGS,alg:SR1}}
    \STATE{Set $\displaystyle 
    Q_0 = \mbox{\tt orth}\left(
    \left[ \; \frac{h_0}{\|\nabla f_0\|} \;,\; \frac{\nabla f_0}{\|\nabla f_0\|} \;,\; Y_0 \;\right]
    \right)$}.
    \COMMENT{See \cref{eq:Qkrep}}
    \STATE{Compute  $P \leftarrow Q_0^\top \,H_0 \,Q_0; \quad
    b \leftarrow  Q_0^\top \,H_0 \,\nabla f_0; \quad
    C \leftarrow Q_0^\top\, H_0^2\, Q_0$.}
    \COMMENT{See \cref{eq:JuliaTRS}}
    \REPEAT[$k=0,1,\ldots$]
    \STATE{Compute $a_k \leftarrow \mbox{{\tt trs\_small}}(P,b,\Delta_k,C)$.}
    \COMMENT{from Julia {\tt TRS} package}
    \STATE{Compute $p_k \leftarrow H_k Q_k a_k$, \quad
    $f_{k+1}\leftarrow f(x_k + p_k)$ and $\rho$ given by \cref{eq:rho}.}
    \STATE{
    Update $\Delta_{k+1}$ according to 
    \cref{eq:TRUpdate}.  
    }
    \IF{ $\rho \le 0$}
    \STATE {Set $x_{k+1} \leftarrow x_k; \quad 
    f_{k+1}\leftarrow f_k; \quad 
    \nabla f_{k+1}\leftarrow\nabla f_k;
    \quad H_{k+1}\leftarrow H_k$.}
    \COMMENT{Reject Step}
    \ELSE
    \STATE{Set $x_{k+1} \leftarrow x_k + p_k$.}
    \COMMENT{Accept Step}
    \STATE{Pick new supplemental directions, 
    $S_{k+1}$, using one of 
    \cref{eq:SUpdate1,eq:SUpdate2,eq:SUpdate3,eq:SUpdate4,eq:SUpdate5,eq:SUpdate6}
    .} 
    \STATE{Compute 
    $\left( \nabla f_{k+1}, h_{k+1},Y_{k+1} \right) 
    \leftarrow \gHS (x_{k+1},S_{k+1})$.}
    \IF{{\em pflag}}
    \STATE{ $H_k \leftarrow \QN(H_k,p_k,\nabla f_{k+1}-\nabla f_k, \delta )$ }
    \ENDIF
    \STATE{Construct 
    $\displaystyle U_{k+1} = \left[\; S_{k+1} \;,\; \frac{\nabla f_{k+1}}{\|\nabla f_{k+1}\|}\; \right]$ and  
    $\displaystyle V_{k+1} = \left[ \; Y_{k+1} \;,\; \frac{h_{k+1}}{\|\nabla f_{k+1}\|}\; \right]$. }
    \STATE{Update
    $H_{k+1} \leftarrow \QN \left( H_{k}, U_{k+1}, V_{k+1}, \delta \right).$}
    \STATE{Set $\displaystyle Q_{k+1} =
    \mbox{\tt orth}\left(
    \left[\;\frac{h_{k+1}}{\|\nabla f_{k+1}\|}\;,\;\frac{\nabla f_{k+1}}{\|\nabla f_{k+1}\|}\;,\; Y_{k+1}\; \right] \right)$}.
    \COMMENT{See \cref{eq:Qkrep}}
    \STATE{Compute \cref{eq:JuliaTRS}\\
    $P \leftarrow Q_{k+1}^\top H_{k+1} Q_{k+1}; \quad
    b \leftarrow  Q_{k+1}^\top H_{k+1} \nabla f_{k+1}; \quad
    C \leftarrow Q_{k+1}^\top H_{k+1}^2 Q_{k+1}$.
      }
    \ENDIF 
    \UNTIL{$||\nabla f_{k+1} || \leq \epsilon $ \COMMENT{Convergence}}
    \RETURN {$x_{k+1}$}
    \end{algorithmic}
\end{algorithm}

\section{Numerical Experiments}
\label{sec:Experiments} 
We test our algorithms on the Rosenbrock function
\begin{align*}
f(x) = \sum_{i=1}^{n} \left[ 
a\left(x_{i+1}-x_i^2\right)^2 + \left(x_i - 1\right)^2\right],
\end{align*}
which is a popular test problem for gradient-based optimization algorithms.
The Julia package used to generate these results in this section is archived \cite{BlockOpt_Code}; a more updated version of the software maybe available at \href{https://github.com/danphenderson/BlockOpt.jl}.
The experiments are initialized with $x_0\in\mathbb{R}^n$ having each component drawn from the uniform distribution on $[-1,1]$.
The global minimum for 
the Rosenbrock function lies in a narrow valley with many saddle points 
(in dimension $60$ 
Kok and Sandrock \cite{RosenSaddles} 
find $53,165$ saddles and predict 
over $145$ million in dimension $100$) 
which makes the minimization challenging for many
algorithms. 
We report results for $n=100$ (a relatively small 
dimensional problem that still illustrates the 
behavior of our algorithms), 
and $a=100$, the standard torture test for 
optimization algorithms. 
We treat gradient evaluations as the primary expense in each 
optimization step and evaluate our algorithms by counting the 
number of $\gHS$ evaluations. Each such evaluation 
involves two 
sequential calls to the underlying AD code $\gAD$ with 
each such call evaluating $w$ simultaneous Hessian samples
in about $2-3$ times an evaluation of $g$.  Given sufficient 
SIMD processors and neglecting linear 
algebra, each algorithmic step takes roughly $5$ times 
the evaluation time for 
a single gradient.  Plots use the number of 
$\gHS$ evaluations. 

\begin{figure}
    \centering
    \includegraphics{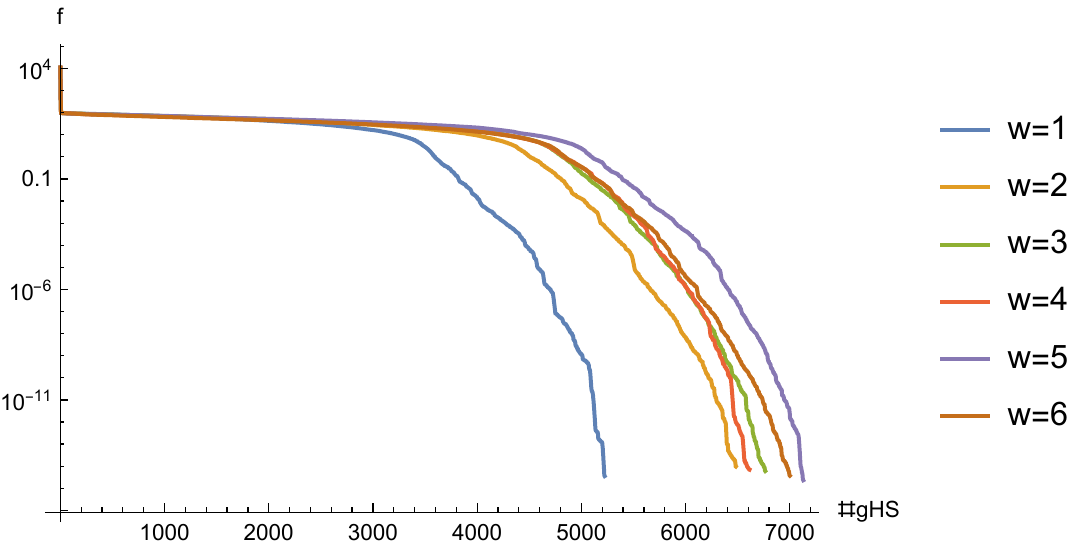}
    \caption{Block PSB with preliminary secant update (i.e., with pflag=1 in \cref{alg:Alg}). Supplemental directions chosen using  \cref{eq:SUpdate4}. Surprisingly, as we increase the sample size, the performance of the algorithm degrades.}
    \label{fig:PSBw}
\end{figure}

\begin{figure}
    \centering
    \includegraphics{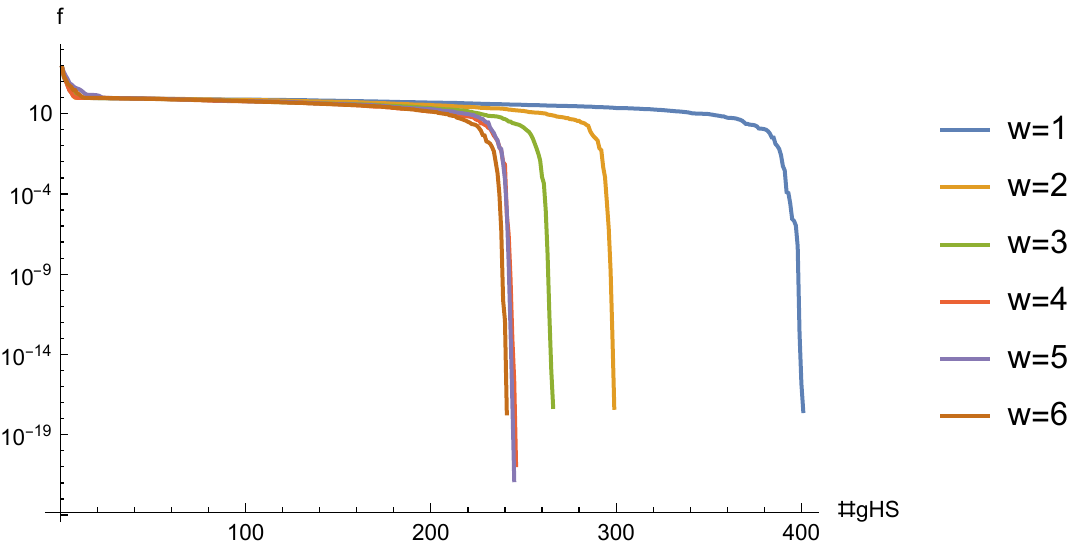}
    \caption{Block SR1 with preliminary secant update  (i.e., with pflag=1 in \cref{alg:Alg}). Supplemental directions chosen using  \cref{eq:SUpdate4}. Comparing \cref{fig:SR1w} and \cref{fig:PSBw}, block SR1 converges much more rapidly than block PSB. Also, the convergence of block SR1 is superior for larger sample sizes $w$, though the improvement appears to decrease as $w$ increases.}
    \label{fig:SR1w}
\end{figure}

\begin{figure}
    \centering
    \includegraphics[width=0.95\textwidth]{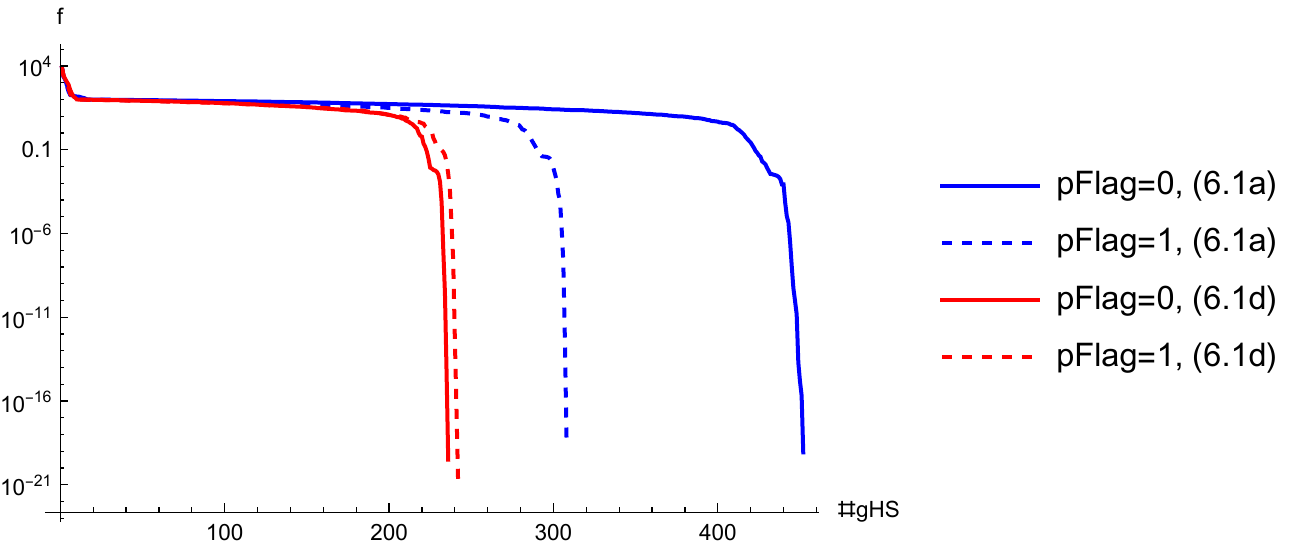}
    \caption{Convergence of SR1 variants with pflag=$\{0,1\}$, and sample directions chosen using $\{\cref{eq:SUpdate1},\cref{eq:SUpdate4}\}$. Sample size is held fixed at $w=4$. The blue curves (sample directions do not include $p_k$, \cref{eq:SUpdate1}) are consistent with the observations in 
    \cite{ByrdSchnableShultzParallelQN88, ByrdSchnableShultzParallelFnEvals88}  that including the approximate secant curvature information in a  preliminary QN update is advantageous. Our framework allows us to 
    incorporate Hessian information in the direction of $p_k$, i.e., \cref{eq:SUpdate4}. The red curves show that using \cref{eq:SUpdate4} to select sample directions is highly beneficial, and eliminates the need to include approximate secant curvature information in the Hessian, which is a serial bottleneck. }
    \label{fig:SR1variants}
\end{figure}

\begin{figure}
    \centering
    \includegraphics{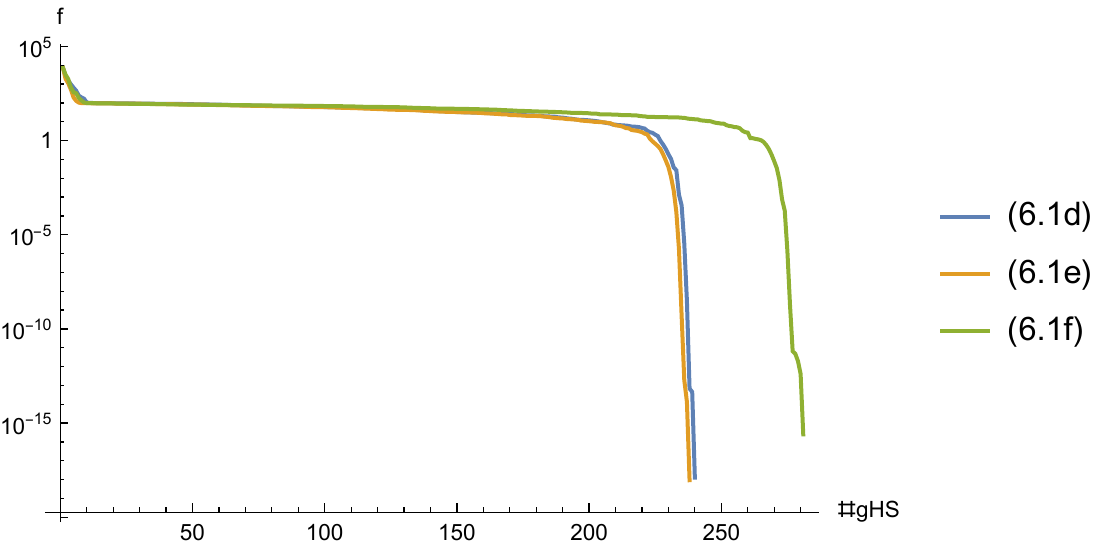}
    \caption{SR1 with with $w=4$, 
    samples from \cref{eq:SUpdate4,eq:SUpdate5,eq:SUpdate6},
    and no preliminary secant update (i.e., with pflag=0 in \cref{alg:Alg}). The simple purely randomized supplemental directions from 
\cref{eq:SUpdate4} gave good results. }
    \label{fig:SampleVariants}
\end{figure}

Representative results for single runs are presented in \cref{fig:PSBw,fig:SR1w,fig:SR1variants,fig:SampleVariants}. At times the 
algorithm converges to the secondary local min described in \cite{RosenSaddles}. These runs are discarded.
\Cref{fig:PSBw} should be compared to \cref{fig:SR1w}
to see the comparatively 
slow convergence of the PSB update. 
The block PSB update is
very conservative in the sense that it gives the smallest
change (in the Frobenius norm) consistent with the new 
Hessian sample which may cause PSB to struggle with
the rapidly changing Rosenbrock Hessian. SR1 is the 
primary focus from here on. As would be expected 
\cref{fig:SR1w} 
shows SR1 converging faster for larger sample sizes 
$w$. The improvement appears to decrease as $w$ increases.

\Cref{fig:SR1variants} compares several SR1 variants with and without a preliminary secant update ({\em pflag} = $\{0,1\}$ in \cref{alg:Alg}) and sample directions chosen from $\{$\cref{eq:SUpdate1,eq:SUpdate4}$\}$. The sample size is held fixed at $w=4$.  The blue curves (sample directions do not include $p_k$, \cref{eq:SUpdate1}) are consistent with the observations in \cite{ByrdSchnableShultzParallelQN88, ByrdSchnableShultzParallelFnEvals88}  that including the approximate secant curvature information in a  preliminary QN update is advantageous. Incorporating $p_k$ in our selection of sample directions,  \cref{eq:SUpdate4} is highly beneficial (red curves), and eliminates the need for the preliminary QN update, which is a serial bottleneck.

Lastly, \cref{fig:SampleVariants} shows the effect 
of varying the sample selection 
strategy between 
\cref{eq:SUpdate4,eq:SUpdate5,eq:SUpdate6} 
The simple purely randomized supplemental directions from 
\cref{eq:SUpdate4} gave good results. The intuition behind \cref{eq:SUpdate5}
was to prevent immediate re-sampling by orthogonalizing against 
the immediate previous directions. It did not lead to significant improvement.
The intuition behind \cref{eq:SUpdate6} was to guide the algorithm to 
accurately resolve eigen-space associated with the larger Hessian eigenvalues.
This variant does appear to resolve these eigenspaces but unfortunately it does
not improve the performance of the algorithm.

\section{Conclusions and Future Work}
\label{sec:Conclusions}
The goal was a straightforward algorithm that would focus on 
potential benefits AD generated Hessian samples in optimization 
algorithms.  The algorithms presented are intended to 
evaluate potential benefits of incorporating block Hessian 
samples in various ways into a fairly standard optimization 
framework.  Practical 
implementations would require limited memory updates
to reduce storage requirements.  Carefully eliminating 
some sampled directions from the trust region sub problem 
and/or exploiting the structure of a limited memory update 
(as shown for LSR1 by Brust et al.~\cite{LSR1TrustProblemBrust2017}) would reduce the computational 
intensity of the trust-region sub-problem solver. We do 
not address these issues in this article and restrict 
attention to significantly fewer than the 
$64$ Hessian samples 
which are feasible on common GPU hardware.

\section{Distribution of Responsibilities}
The article implements matrix 
approximation ideas from 
Dr.~Azzam's thesis in optimization.
Dr.~Struthers designed the algorithm and 
drafted the article with significant input 
from Drs.~Ong and Azzam.  Mr.~Henderson created 
the Julia test problems and code, and
generated numerical results. 
All authors made significant editorial contributions.

\bibliographystyle{alphadin}
\bibliography{references}
\end{document}